\documentclass[11pt]{article}
\usepackage{layout}
\usepackage{amssymb}
\usepackage{float}

\usepackage{graphicx}
\usepackage{epsfig}

\usepackage{color}
\usepackage{latexsym}

\usepackage{amssymb}
\usepackage{algorithmicx}
\usepackage[ruled]{algorithm}
\usepackage{algpseudocode}
\usepackage{algpascal}
\usepackage{algc}

\alglanguage{pseudocode}
\usepackage{graphicx}
\usepackage{amsmath,epsfig,graphics}

\newtheorem{theorem}{Theorem}[section]

\newtheorem{definition}[theorem]{Definition}

\newtheorem{proposition}[theorem]{Proposition}
\newtheorem{remark}[theorem]{Remark}

\newcommand{\Frac}[2] {\frac{\textstyle #1} {\textstyle #2}}

\begin{document}

\title{Parallel Matrix Function Evaluation via Initial value ODE modelling}
\maketitle
\centerline{\scshape Jean-Paul Chehab$^1$ and Madalina Petcu$^{2,3}$ 
}
\medskip
{
\centerline{$^{1}${\footnotesize Laboratoire Amienois de Math\'ematiques Fondamentales et Appliqu\'ees (LAMFA), {\small UMR} 7352}}
  \centerline{{\footnotesize Universit\'e de Picardie Jules Verne, 33 rue Saint Leu, 80039 Amiens France} }
\centerline{$^{2}${\footnotesize UMR 7348   CNRS, Laboratoire de Math\'ematiques et Applications, Universit\'e de Poitiers }} 
 \centerline{{\footnotesize Boulevard Marie et Pierre Curie, T\'el\'eport 2 }}
 \centerline{ {\footnotesize BP30179, 86962 Futuroscope Chasseneuil Cedex, France}}
\centerline{$^{3}${\footnotesize Institute of Mathematics of the Romanian Academy, Bucharest, Romania} }

%
%
\begin{abstract}
The purpose of this article is to propose ODE based approaches for the numerical evaluation of matrix functions $f(A)$, a question of major interest in the numerical linear algebra. To this end, we model $f(A)$ as the solution at a finite time $T$ of a time dependent equation. We use parallel algorithms, such as the parareal method, on the time interval $[0, T]$ in order to solve the evolution equation obtained. When $f(A)$ is reached as a stable steady state,  it  can be computed  by combining parareal algorithms and optimal control techniques. Numerical illustrations are given.
\end{abstract}
{\small 
}\\
\section{Introduction}
The efficient numerical computation of matrix functions, as well as the solution of matrix polynomial equations (such as Riccati's), is one of the actual important topic in numerical linear algebra: this kind of problem arises in a number of situations such as for the approximation of nonlocal operators or of infinitesimal generators as well as in computational control, we refer to \cite{AbreuRaydan,BenziRazouk} for the computation of $p$th-roots of matrices, to  \cite{MonsalvesRaydan1} for the solution of rational matrix equations 
and to the book of Nick Higham  \cite{HighamBook} for a general presentation.\\
Let $A$ be a matrix, to evaluate $f(A)$ or $f(A)b$ for a given vector $b$, number of strategies have been developped: they often use an application of approximation theory techniques to a complex representation of $f$, defining rational approximation methods, let us cite the recent works of Fommer {\it et al} \cite{Fommer} on rational Krylov methods and the references therein. \\
The use of efficients parallel computations of parabolic PDEs have been proposed in e.g. in \cite{GallopoulosSaad} where the generator $\exp(tA)$ is decomposed as a sum of independent polar terms, and more recently, by Gander and Guettel \cite{GanderGuettelParaexp} with the Paraexp algorithm combining Duhamel's principle and additive identities. The evaluation of $\exp(A)$ is an old problem, one can find a nice survey of methods in 
\cite{MolerExponential}. However, in number of situations, the function $f$ is not known so the tools of approximation theory can not be applied. \\

In a more general way, differential equations can be used to identify, then to compute, the matrix function $f(A)$, as a state of the solution; it can be a solution at finite time as well as a stable steady state, when good stability properties are present.\\

To compute $f(A)$ as the solution of an ODE at finite time $T$,  $T=1$ for simplicity, an homotopy method can be displayed as follows: let us consider the matrix differential equation
 \begin{equation} \label{eq1}\begin{split}
 & \Frac{d X}{dt}=\mathcal{F}(X(t)), t \in (0,1),\\
 & X(0)=X_0.
 \end{split}\end{equation}
 We look for $\mathcal{F}$ such that to have $X(1)=f(A)$.  We introduce the homotopy path
 $\mathcal{A}(t)=X_0+t(A-X_0),$
  we set $X(t)=f(\mathcal{A}(t))$; $X_0$ is chosen in such a way $f(X_0)$ is easy to compute.  We have:
  \begin{equation}\label{eqX}\begin{split}
 &\Frac{d X}{dt}= \Frac{d{\mathcal A}}{dt}f'(f^{-1}(X(t))), t \in (0,1),\\
 & X(0)=X_0.
\end{split}\end{equation}
A numerical approximation of $f(A)=X(1)$ can be obtained by using any time marching scheme. For example, taking
$X_0=I$, where $I$ is the identity matrix and ${\mathcal F}(X)=X(Id-A)X$ we have $X(1)=A^{-1}$, so Forward Euler method as well as Runge Kutta method can be applied to build inverse preconditioners of $A$, see \cite{jpc}.\\

However, in a number of situations, it can be difficult to model $f(A)$ as the solution of the ODE at finite time and it is interesting to link $f(A)$ to a stable steady state, see \cite{jpc,DuboisSaidi}. Defining $X=f(A)$ as an asymptotically stable root of ${\mathcal F}$, say ${\mathcal F}(X)=0 \iff X=f(A)$ and assuming that the eigenvalues of the differential of ${\mathcal F}$ at $f(A)$, $D{\mathcal F}(A)$, are of strictly negative real part, we can  consider the differential equation
\begin{equation}\label{stat}
\begin{split}
 & \Frac{d X}{dt}={\mathcal F}(X(t)), t >0,\\
 & X(0)=X_0.
 \end{split}\end{equation}

This modeling allows the computation of sparse approximations to $f(A)$ by defining an associate sparse matrix flow,  see \cite{jpc}.\\

Hence, if a differential system possesses ${\mathcal F}(A)$ as an asymptotical steady state, it is possible to compute numerically ${\mathcal F}(A)$ by applying an explicit time marching scheme. The efficiency of the method depends both on the dynamical properties of the differential system and on the numerical scheme in time.\\

Parallel methods in time have been introduced and developed as a mean to solve time dependent problems using parallel computing, they apply to the computation of the solution of an ODE at a given finite time. These methods where
especially devoted to the computation of $f(A)b$, where $b$ is a given vector in $\mathbb{R}^n$
Let us cite the recent Paraexp and Rational Krylov methods \cite{GanderGuettelParaexp,GuettelRational}, to name but a few. However, it must be noticed that these methods are direct but not well adapted to the computation of $f(A)$, furthermore they assume to have a knowledge of some rational expansions of $f(z)$, $z \in \mathbb{C}$ which is not always the case. The Parareal Method (PM)  is iterative and it can be seen as a multi-steps shooting scheme, where each step can be solved in parallel. Firstly designed for parabolic-like problems,  PM have been successfully adapted to second order evolutive PDE, see  \cite{farhat-06, GanderVdW,GanderPetcu2} but few attention was given for their application in numerical linear algebra and the main purpose of this article is to show that interesting issues can be derived for evaluating $f(A)$.\\

In this article, we consider a framework in which the parallel computation of functions of matrices can be applied, we focus on PM for which we propose  adaptations and implementations  for computing numerical approximations of matrix functions. 
The article is organized as follows: in Section \ref{classical_parareal} we recall some parallel in time methods including the classical parareal algorithm which we apply to solve equation \eqref{eq1} in order to compute the matrix function $f(A)$ of a matrix $A$. We numerically illustrate the algorithm by approximating the inverse and the exponential of a matrix.
In Section \ref{modified_parareal}, we propose a version of the modified parareal algorithm seen as a multiple shooting method as described in \cite{GanderPetcu2}, adapted here to the matricial computations. Using this algorithm, we compute the cosine of a matrix and we compare the results obtained with this algorithm to the case when the cosine is obtained using the classical parareal algorithm. We can notice that the modified parareal algorithm produces an acceleration for the convergence. 
In Section \ref{acceleration} we consider the case when the matrix function is found as the stable steady state of an ordinary differential equation. We propose some methods allowing to accelerate the convergence in time, in order to efficiently compute the steady state.
Section \ref{control} presents a particular acceleration procedure in order to converge to the steady state. The idea is to obtain the matrix function as a solution of an optimal control problem and to apply a parareal in time method in order to solve the control problem, as proposed by \cite{maday-turinici}.

All the numerical computations presented here have been made with Matlab.

\section{Parallel algorithms applied to the computation of matrix functions}\label{classical_parareal}
Before focusing on the application of the parareal method to the numerical evaluation of functions of matrices, we recall briefly some special parallel algorithms in time for solving numerically the equation:
\begin{equation}\label{eq1}\begin{split}
& \Frac{du}{dt}=Au+g,\\
& u(0)=u_0,
\end{split}\end{equation}
and where one needs to be able to evaluate $\exp(tA)b$ for a certain given vector $b$.

A first idea in solving \eqref{eq1} consists in identifying a generator and to use a proper expansion of the underlying function. Namely, following \cite{GallopoulosSaad} we need to approximate the exact solution of \eqref{eq1} which reads as:
\begin{eqnarray}
u(t)=u_0+A^{-1}(\exp(tA)- Id)(Au_0+g).
\end{eqnarray}
Introducing $\phi(t)=\Frac{e^z-1}{z}$, we have $u(t)=u_0+t\phi(tA)(Au_0+g)$. At this point, we consider the polar decomposition of $\phi(s)$ on the spectral interval of $t A$ as
\begin{eqnarray}
\phi(z)\simeq r(z)=\displaystyle{\sum_{j=1}^p\Frac{\omega_j}{s_j-z}}
\end{eqnarray}
Hence, the vector $u$ at time $t$ is computed as:
\begin{equation}
u(t)\simeq u_0 + \displaystyle{\sum_{j=1}^p} t \omega_j (s_j Id - tA)^{-1} (Au_0 + g).
\end{equation}

This last sum can be built in parallel on $p$ processor by solving $p$ independent linear systems
$(s_jId-tA)v_i=A u_0 + g$, and taking then the linear weighted combination $\displaystyle{\sum_{j=1}^p\omega_jv_j}$.\\

Another approach to solve \eqref{eq1} was proposed by Gander and Guettel: the paraexp method. It consists in applying the Duhamel principle writing $u=v+w$ where
\begin{equation}\begin{split}
&\Frac{dv}{dt}=Au,\\
& u(0)=u_0,
\end{split}\end{equation}
and
\begin{equation}\begin{split}
& \Frac{dw}{dt}=Aw+g,\\
& w(0)=0.
\end{split}\end{equation}
Then, decomposing the time interval as $[0,T]=\cup_{j=1}^p[T_{j-1},T_j]$ with $0=T_0<T_1<T_2 <\cdots <T_p=T$, Gander and Guettel define the sequences $v_j$ and $w_j$ as
\begin{itemize}
\item[{\bf Step 1}] For $j=1,\cdots p$ solve numerically with serial integrator the Cauchy problem
$$
v'_j(t)=Av_j(t) +g(t), t\in [T_{j-1},T_j], v_j(T_{j-1})=0
$$
\item[{\bf Step 2}] For $j=1,\cdots p$ solve numerically with an acurate (nearly exact) propagator the Cauchy problem
$$
w'_j(t)=Aw_j(t), \in [T_{j-1},T], w_j(T_{j-1})=v_{j-1}(T_{j-1})
$$
\end{itemize}
Finally the solution is obtained as $u(t)=v_k(t)+\displaystyle{\sum_{j=1}^kw_j(t)}, \mbox{ for } t \in [T_{k-1},T_k]$, see \cite{GanderGuettelParaexp} for more details.\\ 

This last method outperforms the parareal method for solving the original Cauchy problem. Anyway, as stated in the  introduction, it can not be applied in all situations, this is why we first recall the principle of the PM and
apply it for computing the matrix $f(A)$, we do not consider here the problem of evaluating $f(A)b$.

We now concentrate on PM. First introduced for the parallel numerical solution of parabolic problems, we show how to adapt them to the computation of function of matrices.
We start by describing the parareal algorithm that we intend to use. The method was first introduced by Lions, Maday and Turinici in \cite{lions-maday-turinici} and the main goal of this method is to propose a time discretization of a differential evolution equation that allows for parallel implementation, in order to accelerate the computation of the solution, compared to the use of a sequential method. In what follows we introduce the parareal algorithm in the form of a multiple shooting method, as it was described by Gander and Vandewalle in \cite{GanderVdW}. The purpose of the method is to compute in parallel a numerical solution for ordinary differential equations of the form:
\begin{equation}\label{evol1}
\begin{split}
& u'(t)=f(u(t)), \> t \in [0, T],\\
& u(0)=u_0.
\end{split}
\end{equation}

The parareal method uses two propagators: an inexpensive coarse propagator ${\bf G}(T_{n+1}, T_n, x)$ which gives a rough approximation to $u(T_n)$, where $u$ is the solution of equation \eqref{evol1} having $u(T_{n-1})=x$ as initial condition, and an expensive fine propagator ${\bf F}(T_{n+1}, T_n, x)$  which gives a more accurate approximation to the same solution $u(T_n)$. Partitioning the time domain $(0,T)$ into $N$-time subdomains $\Omega_n=(T_n, T_{n+1})$, the algorithm works as follows:

\begin{center}
\begin{minipage}[H]{13.5cm}
  \begin{algorithm}[H]
    \caption{The classical parareal algorithm}\label{cl_parareal}
    \begin{algorithmic}[1]
    \State Step $0$: The algorithm starts with an initial approximation $U^0_n$, $n=0, \dots, N$, which can be found for example by using the coarse propagator:
	\begin{equation}\label{eqStep0}
	U^0_{n+1}={\bf G}(T_{n+1}, T_n, U^0_n), \quad U^0_0=u_0.
	\end{equation}
\State 
    \ Step $k+1$:
	\begin{itemize}
	\item Compute in parallel ${\bf F}(T_{n+1}, T_n, U^k_n)$, for $n=0, \dots, N-1$.
	\item Perform a correction step, using both the coarse and the fine propagators:
	\begin{equation}\label{eqStepk+1}
	U^{k+1}_{n+1}={\bf G}(T_{n+1}, T_n, U^{k+1}_n)+ {\bf F}(T_{n+1}, T_n, U^k_n)-{\bf G}(T_{n+1}, T_n, U^k_n).
	\end{equation}
	\end{itemize}
    \end{algorithmic}
    \end{algorithm}
\end{minipage}
\end{center}

\begin{remark}
	We note that the initial step \eqref{eqStep0} of the algorithm is sequential but not expensive, since it uses only a coarse propagator. At Step $k+1$, correction \eqref{eqStepk+1} provides a more accurate solution by using the fine propagator. The significant advantage of this method, compared to a sequential method using the fine propagator, is that all the expensive computations ${\bf F}(T_{n+1}, T_n, U^k_n)$ are computed in parallel and thus correction \eqref{eqStepk+1} provides a solution as accurate as the one obtained by the sequential use of the fine propagator $F$ but at a cost comparable to a sequential computation using the coarse propagator ${\bf G}$.  
\end{remark}

In order to be able to apply this method to the computation of a matrix function $f(A)$, with $A$ a matrix, we need to be able to see $f(A)$ as the solution at a certain time $T$ of an ordinary differential equation. Thus, we introduce the following matrix function:
\begin{equation}
\label{eqAt}
\mathcal{A}(t)=X_0+t(A-X_0),
\end{equation}
with $X_0$ which can be chosen later and we define $X(t)=f(\mathcal{A}(t))$. The matrix function $X$ satisfies the following ordinary differential equation \eqref{eqX}.
We can note that $f(A)=X(1)$ and thus the numerical computation of $f(A)$ reduces to finding a good approximation for the solution of \eqref{eqX} at time $T=1$.\\

\noindent{\bf Example 1:} The inverse $A^{-1}$ of a matrix $A$

Let $A$ be a matrix belonging to $\mathcal{M}_n(\textbf{R})$ and $I$ the identity matrix in $\mathcal{M}_n(\textbf{R})$. We define $\mathcal{A}(t)=I+ t(A-I)$ which corresponds to taking $X_0=I$ in \eqref{eqAt}. We assume that $\mathcal{A}(t)$ is an invertible matrix for all $t$ in $[0,1]$. Then, the matrix $Q(t)={\mathcal{A}(t)}^{-1}$ satisfies the following ordinary differential equation:
\begin{equation}\label{eqQ}
\begin{split}
& \dfrac{{\rm d}Q}{{\rm d} t}= -Q(t)(A-I)Q(t), \quad t \in I,\\
& Q(0)=I,
\end{split}
\end{equation}
and $A^{-1}=Q(1)$. Our purpose now is to find an approximate solution for the solution $Q$ of \eqref{eqQ} at time $T=1$. 

The idea of finding the inverse of a matrix as the solution at a certain time $T$ of an ordinary differential equation was developped by Chehab in \cite{jpc}. Here the difference is that we solve \eqref{eqQ} by the parareal algorithm described above.

In Figure \ref{inverse} we numerically illustrate the computation  of the inverse of the matrix, using the parareal method. We compare the result obtained by applying to \eqref{eqQ} the classical parareal algorithm to the inverse of the matrix and respectively to the solution obtained by using sequentially a fine propagator.

For the numerical simulations, we first divided the time interval $[0,1]$ into $N=25$ intervals $[T_n, T_{n+1}]$, obtaining thus a coarse grid. Each interval was also divided into $J=200$ subintervals, obtaining a fine grid. For the coarse  and fine propagators, we used an Euler explicit method.  The norm we considered in order to estimate the error at each iteration is ${\rm norm}(\cdot ) =\max_{i, j}({\rm abs} (\cdot) )/ \max_{i,j}({ \rm abs} (A))$, where $A \in \mathcal{M}_n(\textbf{R})$ is the matrix we invert and $i, j=1, \dots, n$. The matrix $A$ that we consider in this example is the $1D$ Laplacian belonging to $\mathcal{M}_{80}(\textbf{R})$. We actually compute the inverse of the matrix $A/2^{10}$; we recover $A^{-1}$ by multiplying the result by $2^{10}$.

 \begin{figure}[h!]
 	\begin{tabular}{cc}
 		 		\includegraphics[width=0.5\textwidth]{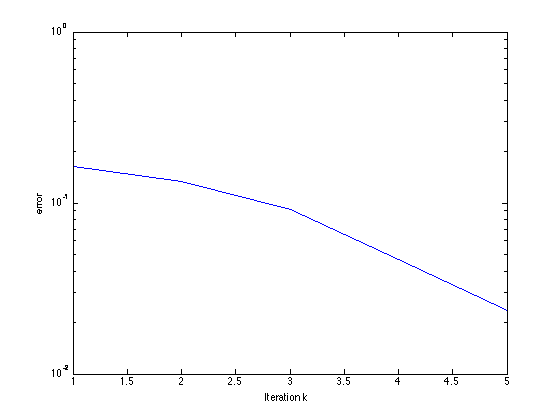}   &
 		\includegraphics[width=0.5\textwidth]{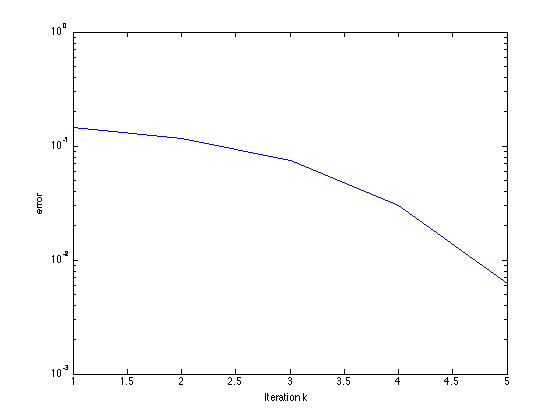}
 	\end{tabular}
 	 	\caption{ Error between the solution given by the parareal algorithm and the inverse of the matrix (left); Error between the solution given by the parareal algorithm and the solution obtained using sequentially a fine propagator (right)}
 	\label{inverse}
 \end{figure}

\noindent{\bf Example 2:} The exponential $\exp(A)$ of a matrix $A$\\
Let $A$ be a matrix belonging to $\mathcal{M}_n(\textbf{R})$ and $I$ the identity matrix in $\mathcal{M}_n(\textbf{R})$. We define $\mathcal{A}(t)=tA$, then the matrix $Q(t)={\exp(\mathcal{A}(t))}$ satisfies the following ordinary differential equation:
\begin{equation}\label{eqQ}
\begin{split}
& \dfrac{{\rm d}Q}{{\rm d} t}= AQ(t), \quad t \in I,\\
& Q(0)=I,
\end{split}
\end{equation}
and $\exp(A)=Q(1)$. As previously, we can find an approximate solution for $Q$ of \eqref{eqQ} at time $T=1$, by using the parareal algorithm. 

In Figure \ref{exp}, we illustrate the computation of $\exp(B)$, where $B=-A/2^{10}$ and $A$ is the Laplacian matrix that we used in the previous example on the inversion of a matrix. We consider the same coarse and fine grids as in the previous example and the propagators are obtained using a Crank-Nicolson method.  The norm we consider in order to estimate the error at each iteration is  ${\rm norm}(\cdot ) =\max_{i, j}({\rm abs} (\cdot) )/ \max_{i,j}({ \rm abs} (\exp(B)))$, where $B \in \mathcal{M}_n(\textbf{R})$ is the matrix for which we compute the exponential and $i, j=1, \dots, n$. We refer the interested reader to \cite{MolerExponential} for a review of methods to compute the exponential of a matrix.

\begin{figure}[h!]
	\begin{tabular}{cc}
		\includegraphics[width=0.5\textwidth]{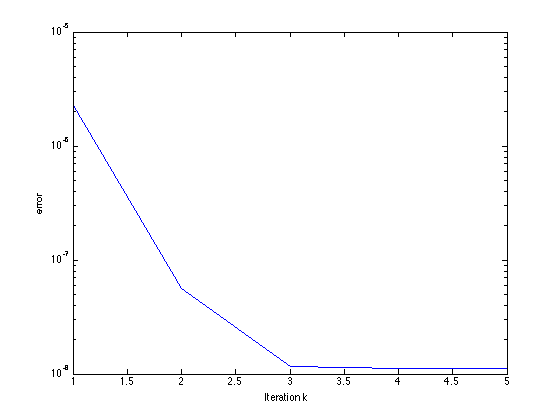}   &
		\includegraphics[width=0.5\textwidth]{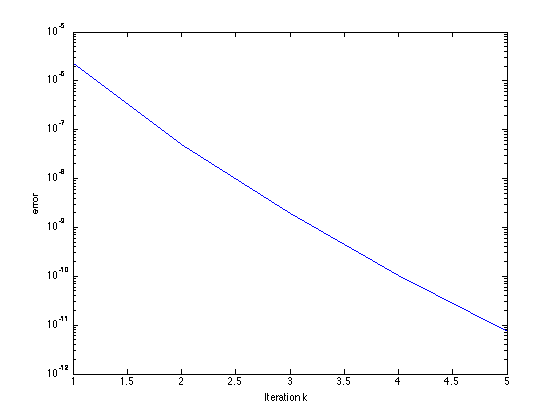}
	\end{tabular}
	\caption{ Error between the solution given by the parareal algorithm and the exponential of the matrix (left); Error between the solution given by the parareal algorithm and the solution obtained using sequentially a fine propagator (right)}
	\label{exp}
\end{figure}

\section{The modified parareal algorithm adapted to the matriceal computations}\label{modified_parareal}

In what follows we adapt the modified parareal algorithm (see \cite{GanderPetcu1}, \cite{GanderPetcu2} for the modified parareal algorithm, see also \cite{farhat-03}, \cite{farhat-06}) to the matriceal case, in order to be able to use it for the computation of matrix functions. 

We start by introducing the Kronecker product and the $\diamond$ product.

 For two matrices $A$ and $B$ belonging to $\mathcal{M}_{n,s}(\mathbf{R})$, we define the inner product 
 $<A, B>_F=tr(B^TA)$, where $tr(M)$ is the trace of a matrix M. The associated norm is the Frobenius norm denoted by $\| \cdot \|_F$. 
 \begin{definition}
 	A system of matrices of $\mathcal{M}_{n,s}(\mathbf{R})$ is said to be $F$-orthonormal if it is orthonormal with respect to the Frobenius inner product $<\cdot, \cdot>_F$.
\end{definition}
 	
 	We also introduce the $\diamond$ product, defined by:
 \begin{definition}
 Let $A=[A_1, A_2, \dots ,A_p]$ and $B=[B_1, B_2, \dots, B_l]$ be matrices belonging respectively to $\mathcal{M}_{n, ps}(\textbf{R})$ and $\mathcal{M}_{n, ls}(\textbf{R})$, where $A_i$, $B_j$ belong to $\mathcal{M}_{n, s}(\textbf{R})$	for all $i=1, \dots, p$ and $j=1, \dots, l$. Then, the $p\times l$ matrix $A^T \diamond B$ is defined by:
$$
A^T\diamond B= \left(\begin{matrix}
 <A_1, B_1>_F &  <A_1, B_2>_F & \cdots & <A_1, B_l>_F\\ 
  <A_2, B_1>_F &  <A_2, B_2>_F & \cdots & <A_2, B_l>_F\\	
   \vdots &  \vdots & \ddots&  \vdots\\
  <A_p, B_1>_F &  <A_p, B_2>_F & \cdots & <A_p, B_l>_F  
 \end{matrix}\right).
$$ 
 \end{definition}
 		
We recall here the following properties for the $\diamond$ product; for more details we refer the reader to \cite{bouyouli}.

\begin{proposition}
	\begin{enumerate}
		\item If $s=1$, then $A^T \diamond B=A^T B$.
		\item The matrix $A=[A_1, A_2, \dots, A_p]$ is $F$-orthonormal if and only if $A^T \diamond A=I_p$.
		\item If $X\in \mathcal{M}_{n,s}(\textbf{R})$, then $X^T \diamond X = \|X\|^2_F$.
	\end{enumerate}
\end{proposition}
In what follows, we also need to introduce the global $QR$-factorization algorithm (see \cite{bouyouli} for more details on the method):\\

\noindent {\large \bf The  global $QR$-factorization algorithm globalQR($\cdot$)}

Let $Z=[Z_1, Z_2, \dots, Z_k]$ be a matrix of $k$ blocks, where $Z_i \in \mathcal{M}_{n,s}(\textbf{R})$, with $i=1, \dots, k$. We search for an $F$-orthonormal family of matrices $\{Q_1, Q_2, \dots, Q_l\}$ such that $Span\{Q_1, Q_2, \dots, Q_l\}=Span\{Z_1, Z_2, \dots, Z_k\}$ and $<Q_i, Q_1>_F=1$ for all $i$, with $<Q_i, Q_j>_F=0$ for all $i$, $j$ such that $i\neq j$.\\

\begin{center}
\begin{minipage}[H]{12cm}
  \begin{algorithm}[H]
    \caption{The globalQR algorithm}\label{global_QR}
    \begin{algorithmic}[1]
    \State Define
    \begin{itemize}
   \item  $R=(r_{i,j})_{i,j=1, \dots, k}=O_k$ 
        \item   $r_{11}=\|Z_1\|_F$ 
        \item $Q_1=Z_1/r_{11}$
           \end{itemize}
    \For{  $i=2, \dots, k$}
    \begin{itemize}
      \item  $Q=Z_i$
	        \For{$j=1, \dots, i-1$}\\
	       $r_{j,i}=<Q, Q_j>_F, \quad Q=Q-r_{j,i}Q_j$
	       \EndFor  $\>j$
	     \item  $r_{ii}=\|Q\|_F$\\
	       If $r_{ii}\neq 0$, then $Q_i=Q/r_{ii}$
    \end{itemize}
     Else $Q_i=Q$
            \EndFor  $\>i$
            \State If $r_{ii}=0$, eliminate $Q_i$ from the family $\{Q_1, \dots, Q_k \}$ and eliminate the $i$-th line and column from $R$
    \end{algorithmic}
    \end{algorithm}
\end{minipage}
\end{center}

The result of the algorithm is $[Q, R]={\rm globalQR(Z)}$.

\begin{proposition}
	\begin{enumerate}
		\item The algorithm $(1)-(4)$ provides the matrix $Q=[Q_1, \dots, Q_k]$ such that
		$$
		Z=[Z_1,\dots, Z_k]=Q(R \otimes I_s),
		$$
		where
		$$
		Q^T\diamond Q=\left( \begin{matrix}
		1 &\\
		& 1 & \\
		& & \ddots \\
		& & & 0 & \\
		& & &  & \ddots \\
		& & & & & 1\\
		\end{matrix} \right),
		$$
		the coefficient $0$ on the diagonal of $Q^T \diamond Q$ corresponding to the $i$-th position on which $r_{ii}=0$.
		\item  The algorithm $(1)-(5)$ provides a basis  $\{Q_1, \dots, Q_l\}$ in $Span\{Z_1, \dots, Z_k\}$.
	\end{enumerate}
\end{proposition}

Using the global QR-factorization algorithm introduced above, we are now able to construct the projection on the space generated by a family of matrices.\\

\noindent {\large \bf The projection on a space generated by a family of matrices $\{Z_1, \dots, Z_k\}$}

Let $\{Z_1, \dots, Z_k\}$ be a family of matrices. By applying the global QR-algorithm described above to the family $\{Z_1, \dots, Z_k\}$, we obtain $\{Q_1, \dots, Q_l\}$ to be a family of F-orthonormal generators of $Span\{Z_1, \dots, Z_k\}$.

We define $Q=[Q_1, \dots, Q_l]$. Then, $Q^T \diamond Q$ gives a vector $(\alpha_1, \dots, \alpha_l)$ which defines the projection $P_Q Y$ of $Y$ on the matrix space $Span\{Z_1, \dots, Z_k\}=Span\{Q_1, \dots, Q_l\}$:
$$
P_QY=\sum_{i=1}^l \alpha_i Q_i,
$$
with $(\alpha_1, \dots, \alpha_l)^T=Q^T \diamond Y$.\\

\noindent{\large \bf The modified parareal algorithm}

Taking into account these tools, we are now able to adapt the modified parareal algorithm described in \cite{GanderPetcu2} to the matrix case, in order to be able to solve in parallel the following evolution equation:
\begin{equation}\label{eq_evol}\begin{split}
& U'(t)=f(U(t)), \quad t \in [0, T],\\
& U(0)=U_0,
\end{split}
\end{equation}
with $U: \textbf{R} \rightarrow \mathcal{M}_{n,p}(\textbf{R})$, $U_0 \in \mathcal{M}_{n,p}(\textbf{R})$ and $f : \mathcal{M}_{n,p}(\textbf{R}) \rightarrow \mathcal{M}_{n,p}(\textbf{R})$. 

We will start by studying the linear homogeneous case $U'(t)=\mathcal{B} U(t)$, with $\mathcal{B} \in \mathcal{M}_{n,n}(\textbf{R})$.\\

\noindent {\bf The linear homogeneous case}

The parareal algorithm that we present here is a multiple shooting method using two propagation operators. As for the classical parareal algorithm, we start by splitting the time interval $[0, T]$ into $N$-time subdomains $\Omega_n=(T_{n-1}, T_n)$, with $n \in \{1, \dots, N\}$ and $T_n - T_{n-1}=\Delta T_n$. As previously, let $G(T_n, T_{n-1}, x)$ be a coarse propagator which gives a rough approximation to $u(T_n)$, where $u$ is the solution of equation \eqref{eq_evol} having $u(T_{n-1})=x$ as initial condition and let $F(T_n, T_{n-1}, x)$ be a more expensive fine propagator which gives a more accurate approximation to the same $u(T_n)$. The coarse propagator can be any sequential method while the fine propagator can be obtained by applying a sequential numerical method on the time interval $\Omega_n$ which was partitioned into $J$ subintervals of size $\Delta t_n$ (obtaining thus a fine grid for the whole time interval $[0, T]$). The algorithm works as follows:
\begin{center}
\begin{minipage}[H]{15.5cm}

  \begin{algorithm}[H]
    \caption{The modified parareal algorithm for linear homogeneous problems}\label{modif_parar_hom}
   
    \begin{algorithmic}[1]
    \State Step 0: The algorithm starts by constructing an initial approximation $U^0_n$, $n \in \{ 1, \dots, N\}$, which can be obtained for example by using the coarse propagator sequentially:
	\begin{equation}
	U^0_{n+1}=G(T_{n+1}, T_n, U^0_n), \quad U^0_0=U_0.
	\end{equation}
	We also define the space $\mathcal{S}^{-1}$ as the empty space.
           \State   Step $k+1$:
	\begin{itemize}
		\item We enhance the subspace $\mathcal{S}^k=Span\{\mathcal{S}^{k-1} \cup \{ U^k_n, \> n=0, \dots, N\}\}$, with $\{U^k_n\}_n$ known from the previous step.
		\item We find a basis in $\mathcal{S}^k$ using the global QR-algorithm. 
		
		In order to do so, we construct 
		$$
		Z^k=[S^{k-1}, U_0^k, \dots, U_{N}^k],
		$$
		where $S^{k-1}$ is the matrix having as columns the basis of $\mathcal{S}^{k-1}$. Applying the global QR-factorization, we find the matrices $Q^k$ and $R^k$ such that $[Q^k, R^k]=globalQR(Z^k)$.
		\begin{remark}
			We note here that $Q^k=[Q^{k-1}, Q^k_{j_1}, \dots, Q^k_{j_l}]$ since the basis in $\mathcal{S}^k$ is obtained by completion of the basis in $\mathcal{S}^{k-1}$.
		\end{remark}
		\item We compute in parallel the evolution of $\{Q^k_{j_i}\}_{j_i}$ through the fine propagator:
		$$
		F(T_{n+1}, T_n, Q^k_{j_i}), \quad j_i \in \{j_1, \dots, j_l\}, \quad n\in \{0, \dots, N-1\}.
		$$
		\begin{remark}
			At this step we only need to compute the evolution of $Q^k_{j_i}$ with $j_i \in \{j_1, \dots, j_l\}$ since the evolution of all the other elements of the basis is known from previous steps.
		\end{remark}
		\item We can now perform the sequential update:
		\begin{equation}\label{eq_update}
		U^{k+1}_{n+1}=F(T_{n+1}, T_n,P^k U^{k+1}_n) + G(T_{n+1}, T_n,(I-P^k) U^{k+1}_n),
		\end{equation}
		with $P^k$ the projection on the space $\mathcal{S}^k$.
		
		In order to compute $F(T_{n+1}, T_n,P^k U^{k+1}_n)$, we first compute 
		$$
		(Q^k)^T \diamond U^{k+1}_n =\left(\begin{matrix}
		\alpha_1\\ \vdots \\ \alpha_s
		\end{matrix}\right),
		$$
		and thus $P^k U^{k+1}_n = \sum_i \alpha_i Q^k_i$. Then, the evolution of $P^k U^{k+1}_n$ through the fine propagator can be found as:
		$$
		F(T_{n+1}, T_n,P^k U^{k+1}_n)=\sum_i \alpha_i F(T_{n+1}, T_n,P^k Q^{k}_i).
		$$
	\end{itemize}       
	\end{algorithmic}
    \end{algorithm}
\end{minipage}
\end{center}

\begin{remark}
	\begin{enumerate}
		\item The only sequential computation in \eqref{eq_update} is ${\bf G}(T_{n+1}, T_n,(I-P^k) U^{k+1}_n)$ which is a cheap computation since it uses only a coarse propagator. No extra evaluation using the fine propagator is required in \eqref{eq_update} since the fine evolution on the space $\mathcal{S}^k$ is already known and has been computed in parallel.
		\item Note that if we formally take $k \to \infty$ in \eqref{eq_update}, then the projection $P^k$ tends to the identity matrix and thus $U^{k+1}_n$ tends to the approximate solution given by the fine propagator. A rigorous explanation for this phenomenon can be found in \cite{GanderPetcu2}.
	\end{enumerate}
\end{remark}

\noindent {\bf The linear inhomoneneous case}\\
The algorithm for the inhomogeneous case follows the same idea as in \cite{GanderPetcu1}, \cite{GanderPetcu2}. Denoting by ${\bf F}$ the matrix that defines the fine propagator, we can notice that ${\bf F}(T_{n+1}, T_n, \alpha X + \beta Y)$can be computed as follows:
\begin{equation}\begin{split}
&{\bf F}(T_{n+1}, T_n, \alpha X + \beta Y)=\alpha {\bf F}^J X + \beta {\bf F}^J Y + {\bf F}(T_{n+1}, T_n, 0)\\
& \> =\alpha ({\bf F}(T_{n+1}, T_n, X)- {\bf F}(T_{n+1}, T_n, 0))+ \beta ({\bf F}(T_{n+1}, T_n, Y)- {\bf F}(T_{n+1}, T_n, 0))\\
&\quad \>+{\bf F}(T_{n+1}, T_n, 0).
\end{split}\end{equation}
Thus, the affine computation for ${\bf F}(T_{n+1}, T_n, P_k U^{k+1}_n)$, and similarly  for ${\bf G}(T_{n+1}, T_n, (I-P_k) U^{k+1}_n)$, requires only an additional computation, which leads to the following algorithm:

\begin{center}
\begin{minipage}[H]{15.5cm}

  \begin{algorithm}[H]
    \caption{The modified parareal algorithm for linear inhomogeneous problems}\label{modif_parar_inhom}
    
    \begin{algorithmic}[1]
    \State Step $-1$: We start the algorithm by computing the quantities ${\bf F}(T_{n+1}, T_n, 0)$ and ${\bf G}(T_{n+1}, T_n, 0)$. These quantities will be useful for the affine computation of $F(T_{n+1}, T_n, P_k U^{k+1}_n)$ and ${\bf G}(T_{n+1}, T_n, (I-P_k) U^{k+1}_n)$ at step $k+1$.
           \State   Step $0$: This step is identical to Step $0$ for the homogeneous case.
        \State Step $k+1$: \begin{itemize}
		\item Enhance the space $\mathcal{S}^k=Span\{\mathcal{S}^{k-1} \cup \{U^k_n, \> n=0, \dots, N \}\}$ as in the homogeneous case and compute the matrix $Q^k=[Q^{k-1}, Q^k_{j_1}, \dots, Q^k_{j_l}]$ in order to have the projection $P^k$. 
		\item Compute in parallel ${\bf F}(T_{n+1}, T_n, Q^k_{j_i})$, for $j_i=j_1, \dots, j_l$.
		\item Update the approximate solution using the formula:
		\begin{equation}\label{update2}
		U^{k+1}_{n+1}={\bf F}(T_{n+1}, T_n, P_k U^{k+1}_n) + {\bf G}(T_{n+1}, T_n,(I- P_k) U^{k+1}_n)- {\bf G}(T_{n+1}, T_n, 0).
			\end{equation}
	\end{itemize}            
	\end{algorithmic}
    \end{algorithm}
\end{minipage}
\end{center}

\begin{remark}
	We note here that ${\bf F}(T_{n+1}, T_n, P_k U^{k+1}_n)$ can be computed at minimal cost, since $P_k U^{k+1}_n$ is a linear combination of approximations for which the evolution is known. In fact, ${\bf F}(T_{n+1}, T_n, P_k U^{k+1}_n)$ is computed as a linear combination of ${\bf F}(T_{n+1}, T_n, 0)$, computed at Step $-1$, and of ${\bf F}(T_{n+1}, T_n, Q^k_j)$, computed in parallel:
	\begin{equation}
	{\bf F}(T_{n+1}, T_n, P_k U^{k+1}_n)=\sum_i \alpha_i \{{\bf F}(T_{n+1}, T_n, Q^k_i)-F(T_{n+1}, T_n, 0)\}+ {\bf F}(T_{n+1}, T_n, 0),
	\end{equation}
	where $P_k U^{k+1}_n = \sum_i \alpha_i Q^k_i$.\\
\end{remark}

\noindent{\large \bf Application of the modified parareal algorithm to the computation of matrix functions}\\

The purpose of this section is to propose an alternative method, using the modified parareal algorithm, to the computation of matrix functions $f(A)$, such as $\sin(A)$, $\cos(A)$, $A^{-1}$, where $A$ is a matrix. As already explained previously, the idea is to see $f(A)$ as the solution at time $T=1$ of a differential equation and to apply the modified parareal algorithm in order to compute the solution of the differential equation.

In order to obtain the differential equation, we can pose $\mathcal{A}(t)=X_0 + t(A-X_0)$, with $X_0$ conveniently chosen, and we define $X(t)=f(\mathcal{A}(t))$. Thus, $X(1)=f(A)$, so we need to compute $X(t)$ at $t=1$. Under optimal conditions, the differential equation satisfied by $X$ is \eqref{eqX}.\\

\noindent{\bf Example 3:} The computation of sine and cosine of a matrix $A$\\
Let $A$ be a given matrix, we write $\mathcal{A}(t)=X_0 + t(A-X_0)$ and we define $X(t)=\sin(\mathcal{A}(t))$, $Y(t)=\cos(\mathcal{A}(t))$. Then, $X$, $Y$ satisfy the following differential equation:
\begin{equation}\begin{split}
& \dfrac{{\rm d}X}{{\rm d}t}=(A-X_0) Y(t),\\
& \dfrac{{\rm d}Y}{{\rm d}t}=-(A-X_0) X(t),
\end{split}
\end{equation}
which can be written in the compact form:
\begin{equation}\label{eqB}
\dfrac{{\rm d}U}{{\rm d}t}=\mathcal{B} U,
\end{equation}
with $U(t)=(X(t), Y(t))$ and $ \mathcal{B}=\left( \begin{matrix}
0 & A-X_0\\
-(A-X_0) & 0
\end{matrix}\right)$;
$X_0$ can be chosen for convenience either $0$ or $I$. 

We can now apply the modified parareal algorithm to the linear homogeneous problem \eqref{eqB} with the initial condition $U(0)=(X(0), Y(0))$. In order to improve the efficiency of the method when we compute $\cos(A)$, we combine the algorithm proposed by Hargreaves and Higham in \cite{HighamSine} and the modified parareal algorithm. Our algorithm reads as follows:\\

\begin{center}
\begin{minipage}[H]{12cm}
  \begin{algorithm}[H]
    \caption{Modified parareal algorithm applied to the computation of  the cosine of a matrix $A$}\label{Alg_cos}
    \begin{algorithmic}[1]
    \State $A$ given in ${\cal M}_n({\mathbb R})$
        \State   Find the smallest non-negative integer $m$ such that 
	$2^{-m} \|A\|_\infty \leq 1$.
        \State Define the matrix $A_0=2^{-m} A$.
             \State Compute $C_0=\cos(A_0)$ using the modified parareal algorithm described previously.
              \State Recover $\cos(A)$ using $m$ times the formula $\cos(2M)=2\cos^2(M)-I$:
         \begin{itemize}
          \For{  $i=0, \dots, m-1$}\\
	 $C_{i+1}=2C_i^2-I$
            \EndFor
        \item    $X=C_m$
            \end{itemize}
            \end{algorithmic}
    \end{algorithm}
\end{minipage}
\end{center}

We illustrate numerically these results. In Figure \ref{cos}, we compare the speed of convergence between the classical parareal algorithm and the modified one, both applied to the computation of the cosine of a matrix. 
For the numerical simulations, we first divided the time interval $[0,1]$ into $N=10$ intervals $[T_n, T_{n+1}]$, obtaining thus a coarse grid. Each interval was also divided into $J=100$ subintervals, obtaining a fine grid. For the coarse  and fine propagators, we used an Euler explicit method. 

We plotted the error measured in the $L^\infty(0,T)$-norm between the fine solution, meaning between an approximation of $\cos(A)$ computed as the solution of the ordinary differential equation \eqref{eqB} when a sequential method using the fine propagator is used, and the approximation produced by the two algorithms at each iteration. 
We can notice that the modified parareal algorithm converges much faster. 

 \begin{figure}[h!]
 	\begin{tabular}{cc}
 		\includegraphics[width=0.5\textwidth]{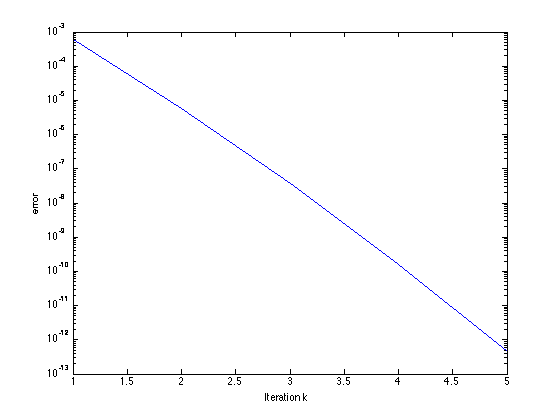}   &
 		\includegraphics[width=0.5\textwidth]{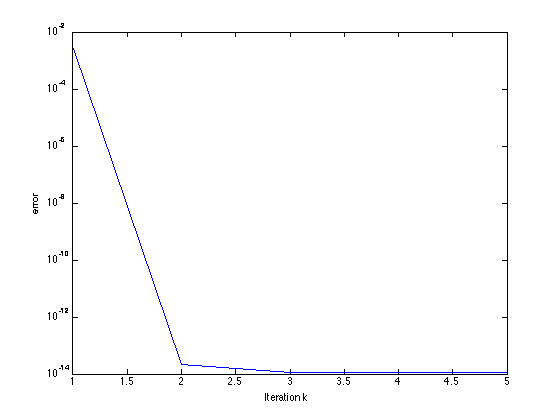}
 	\end{tabular}
 	\caption{ Error for the classical and the
 		modified parareal algorithm, applied to the computation of the cosine of a matrix}
 	\label{cos}
 \end{figure}
\section{Acceleration of convergence in time}\label{acceleration}
\subsection{Description of the method}
In a number of cases, the matrix function evaluated at $A$, $f(A)$, can be defined implicitely as the zero of a proper functional $\mathcal{F}$, i.e., $\mathcal{F}(X,A)=0\iff X=f(A)$; it can correspond to a (stable) steady state of an ODE: 
 \begin{equation}\begin{split}
& \Frac{d X}{dt}= \mathcal{F}(X,A), t>0,\\
 & X(0)=X_0.
 \end{split}\end{equation}
So $f(A)$ can be approached  numerically by applying an explicit time marching scheme on a sufficiently large time interval $[0,T]$. This is the case, e.g., with $\mathcal{F}(X,A)=I-AX$ when considering the computation of $f(A)=A^{-1}$, see \cite{jpc}.

Of course, the efficiency of the numerical method is related to the speed of convergence in time to the steady state. A first idea is to speed it up by adding a time dependent term $V(t)$ as
\begin{equation}\label{ODEaccconv}\begin{split}
 & \Frac{d Y}{dt}= \mathcal{F}(Y,A)+V(t), t>0,\\
 & Y(0)=X_0,
\end{split}\end{equation}
in such a way
\begin{eqnarray}\label{AccConvRel}
 \displaystyle{\lim_{t\rightarrow +\infty}\Frac{\parallel {\mathcal F}(Y(t),A)\parallel}{\parallel {\mathcal F}(X(t),A)\parallel}}=0;
 \end{eqnarray}
 this is in fact a kind of acceleration of convergence in time. An important feature is that superlinearly convergent algorithms can be then obtained by applying appropriate time marching schemes, as shown ahead in this section  with the numerical illustrations.
\begin{remark}
When choosing ${\mathcal F}(Y,A)=f-AY$, the discretization of the above systems by a variable time step-size Backward Euler method gives rise to the iterations
\begin{eqnarray}
X^{(k+1)}=X^{(k)}+\alpha_k\left(f-AX^{(k)}\right),
\end{eqnarray} 
 and respectively
  \begin{eqnarray}
  Y^{(k+1)}=Y^{(k)}+\alpha_k\left(f-AY^{(k)} + V^{(k)}\right).
  \end{eqnarray} 
So, in order to have a solution such that $\displaystyle{\lim_{k\rightarrow +\infty}\Frac{\parallel {\mathcal F}(Y^{(k)},A)\parallel}{\parallel {\mathcal F}(X^{(k)},A)\parallel}}=0$ we need to take $\alpha_k=\Frac{<r^k,r^k>}{<r^k,Ar^k>}$, with $r^k=f-AY^{(k)}$  and $V^{(k)}$ the sequence of the conjugate direction of the conjugate gradient method (CGM), thanks to its convergence in a finite number of steps.
More generally, the solution of linear systems can be modeled by dynamical systems and numerical algorithms can be derived by numerical time integration, see \cite{cl}. 
\end{remark}
We go back to the continuous acceleration of convergence condition (\ref{AccConvRel}). Let us consider for the maximum of simplicity the scalar differential equation
\begin{equation}\begin{split}
& \Frac{dx}{dt}=b-ax,\\
& x(0)=x_0,
\end{split}\end{equation}
where $a>0$.  For any initial datum $x_0$, the solution $x(t)$ converges to $x^*=\Frac{b}{a}$ and we have
$x(t)=e^{-at}x_0+b\Frac{1-e^{-ta}}{a}$. Now we want to tune the equation in such a way the convergence to $x^*$ is faster. This means to determine $u(t)$ such that the solution $y$ of the ODE
\begin{equation}\begin{split}
& \Frac{dy}{dt}=b-ay+u(t),\\
& y(0)=x_0,
\end{split}\end{equation}
satisfies
$\displaystyle{\lim_{t \rightarrow +\infty}\Frac{\mid ay(t)-b \mid}{\mid ax(t)-b\mid}}=0$. A trivial computation gives
$$
\Frac{ay(t)-b}{ax(t)-b}=1+a\Frac{\displaystyle{\int_0^t e^{as}u(s)ds}}{ax_0-b}.
$$
Hence setting $v(s)=e^{as}u(s)$, the last condition reads
$$
\displaystyle{\lim_{t \rightarrow +\infty}\displaystyle{\int_0^t v(s)ds}=\Frac{b-ax_0}{a}}.
$$
Of course the choice of $v(t)$ is not unique. In particular the speed of convergence of $y$ to $x^*$ depends on the speed of convergence of the integral to the limit $\Frac{b-ax_0}{a}$; compactly supported functions $v$, which allow to reach the limit in finite time,  can be good candidates.\\

This little computation can be applied in higher dimension considering a diagonalizable matrix $A$ as
\begin{equation}\begin{split}
& \Frac{dX}{dt}=b-AX,\\
& X(0)=X_0,
\end{split}\end{equation}
and
\begin{equation}\begin{split}
& \Frac{dY}{dt}=b-AY+u(t),\\
& Y(0)=X_0.
\end{split}\end{equation}
If $D$ denotes the diagonal matrix containing the eigenvalues of  $A$ and $P$ the associated passage matrix, we obtain as a condition of acceleration of convergence
$$
\displaystyle{\lim_{t \rightarrow +\infty}\displaystyle{\int_0^t e^{\lambda_i s}{\hat u}_i(s)ds}=\Frac{{\hat b}_i-\lambda_i{\hat X}_{0,i}}{\lambda_i}},
$$
which reads as (\ref{AccConvRel}) choosing an appropriate norm. 
We have used here the notation ${\hat X}=P^{-1}X$ and $D=P^{-1}AP$, $D=diag(\lambda_1,\cdots, \lambda_n)$, so we recover the condition $\displaystyle{\lim_{t \rightarrow +\infty}\displaystyle{\int_0^t v(s)ds}=X^*-X_0}$, with
$v(s)=e^{sA}u(s)$.\\

The above computation can be extended in the nonlinear case as following:
\begin{equation}\begin{split}
&\Frac{dX}{dt}={\cal F}(X),\\
& X(0)=X_0.
\end{split}\end{equation}
We assume that ${\cal F}$ is differentiable at its root $X^*$ and that $D{\cal F}(X^*)$ is diagonalisable and that its eigenvalues of  are of strictly negative real part. We assume in addition that $X_0$ is chosen in a suitable neigborhood of $X^*$ and consider the two linear problems
\begin{equation}\begin{split}
& \Frac{d\delta X}{dt}=D{\cal F}(X^*)\delta X,\\
& \delta X(0)=X_0-X^*,
\end{split}\end{equation}
and
\begin{equation}\begin{split}
& \Frac{d\delta Y}{dt}=D{\cal F}(X^*)\delta Y+U(t),\\
& \delta Y(0)=X_0-X^*.
\end{split}\end{equation}
We deduce from above a way to compute $U(t)$ as
$$
\displaystyle{\lim_{t \rightarrow +\infty}\displaystyle{\int_0^t e^{sD{\cal F}(X^*)}U(s)ds}}=X^*-X_0.
$$
The function $Y(t)=X^0+\delta Y(t)$ is expected to converge to $X^*$ faster than $X(t)$, at least locally. Of course $X^*$ is unknown but it can be estimated during the numerical computation.
\begin{remark}
We can write the acceleration convergence condition as
$$
\displaystyle{\int_0^t e^{sM}U(s)ds}=X^*-X_0 + \xi (t),
$$
with $\displaystyle{\lim_{t \rightarrow +\infty}\xi(t)}=0$; here $M$ is a given matrix, it can be $A$ as well as
$D{\cal F}(X^*)$, depending on the context.
We must impose $\xi(0)=-(X^*-X_0)$ and write $\xi(t)=\displaystyle{\int_0^t \eta(s)ds}-(X^*-X_0)$. Finally, we have
$$
\displaystyle{\int_0^t e^{sM}U(s) -\eta(s) ds}=0 \ \forall t,
$$
so 
$$
U(s)=e^{-sM}\eta(s).
$$
This very simple computation shows how to rely the speed of convergence given by $\xi(t)$ then by $\eta(t)$ to
the accelerator $U(s)$. 
\end{remark}
\subsection{Practical implementation and illustration}
\noindent {\bf Example 4: } The linear vectorial case\\
We first consider the linear problem in ${\mathbb R}^n$ and to illustrate the convergence acceleration, we choose as
control function $u(t)=\chi_{[0,1]}(t)e^{-tA}(X^*-X_0)$ in such a way $\displaystyle{\lim_{t \rightarrow +\infty}\displaystyle{\int_0^t e^{sA}u(s)ds}=X^*-X_0}$. For this choice of $u$ 
we can compute the speed of convergence $\xi(t)$ of $Y(t)$ to $X^*$. Indeed, we have $\eta(t)=\chi_{[0,1]}(t)(X^*-X_0)$. Hence  $\xi(t)=min((t-1),0)(X^*-X_0)$.\\
Of course, the definition of $u$ needs the solution $X^*$ to be known and we will rather use an approximation of this function by replacing $X^*$ with an approximation 
${\tilde X}$ obtained, e.g. by solving a system $M{\tilde X} =b$, where $M$ is a preconditioner of $A$.
\begin{center}
\begin{minipage}[H]{12cm}
  \begin{algorithm}[H]
    \caption{Simple Gradient Acceleration}\label{SGA}
    \begin{algorithmic}[1]
    \State $X^{(0)}=X_0$ given in ${\mathbb R}^n$
        \State  ${\tilde X} \approx X^*$ is approached by solving $M{\tilde X} =b$
        \State $e={\tilde X}-X^{(0)}$ and $u^{(0)}=e$
            \For{$k=0,1, \cdots$ until convergence}
             \State Set $X^{(k+1)}=X^{(k)} +\Delta t(b-AX^{(k)} +u^{(k)})$
              \State Set   $u^{(k+1)}=(Id-\Delta t A)u^{(k)}$ 
            \EndFor
    \end{algorithmic}
    \end{algorithm}
\end{minipage}
\end{center}
Numerically, the limit is reached for a small time $t^*$, that can be determinated. It can be a starting point for implementing the PM method: a coarse solution is first computed to determine $t^*$ and ${\tilde X}$. Then the PM process can be displayed combining fine and coarse propagators.\\
Similarly, we can consider the acceleration of the steepest descend method by taking a sequence of variable time steps as follows:
\begin{center}
\begin{minipage}[H]{12cm}
  \begin{algorithm}[H]
    \caption{Steepest Descent Acceleration}\label{SDA}
    \begin{algorithmic}[1]
    \State $X^{(0)}=X_0$ given in ${\mathbb R}^n$
        \State  ${\tilde X} \approx X^*$ is approached by solving $M{\tilde X} =b$
        \State $e={\tilde X}-X^{(0)}$ and $u^{(0)}=e$
        \State $r^{(0)}=b-AX^{(0)}$.
            \For{$k=0,1, \cdots$ until convergence}
            \State $\Delta t=\Frac{<r^{(k)},r^{(k)}>}{<Ar^{(k)},r^{(k)}>}$
             \State Set $X^{(k+1)}=X^{(k)} +\Delta t (r^{(k)} +u^{(k)})$
              \State Set   $u^{(k+1)}=(Id-\Delta t A)u^{(k)}$ 
              \State Set $r^{(k+1)}=r^{(k)}-\Delta t A(r^{(k)} + u^{(k)})$
            \EndFor
    \end{algorithmic}
    \end{algorithm}
\end{minipage}
\end{center}
Below, on Figures \ref{ACC_SYSLIN_GRAD1} and \ref{ACC_SYSLIN_GRAD2} we reported the history of the ratio 
$\Frac{\parallel AY(t)-b\parallel }{\parallel AX(t)-b\parallel}$ versus the time, when using respectively the simple gradient method  and the steepest descent method. In both cases, the acceleration procedure allows a much faster convergence to the limit. It has to be pointed out that the underlying idea is based only on dynamical system arguments and not on linear algebra ones.
\begin{figure}[!h]
\begin{center}
\includegraphics[width=7.5cm, height=6.5cm]{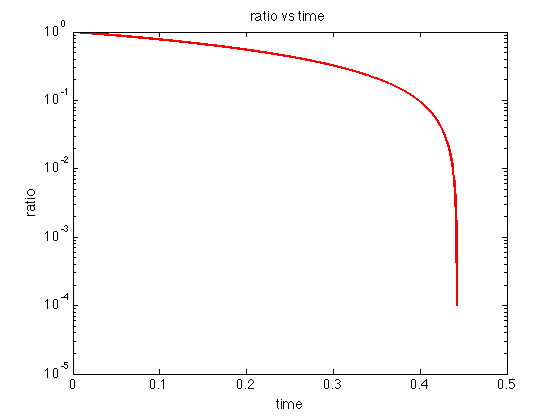}
\includegraphics[width=7.5cm, height=6.5cm]{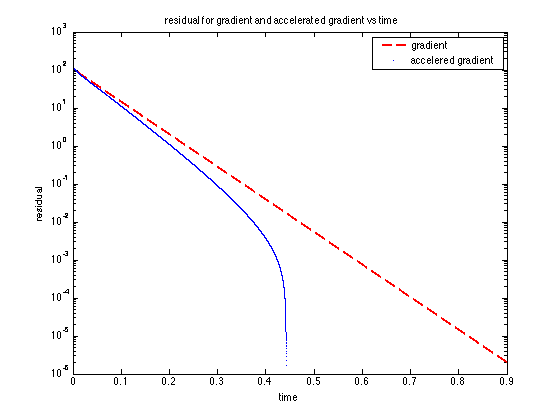}\\
\caption{Acceleration of the convergence in time for the Heat equation $n=127$, $\Delta t=\Frac{h^2}{4}$, $X_0=0$, $b=1$, $M$ is built as an incomplete $LU$ factorization of $A$}
\label{ACC_SYSLIN_GRAD1}
\end{center}
\end{figure}
\begin{figure}[!h]
\begin{center}
\includegraphics[width=7.5cm, height=6.5cm]{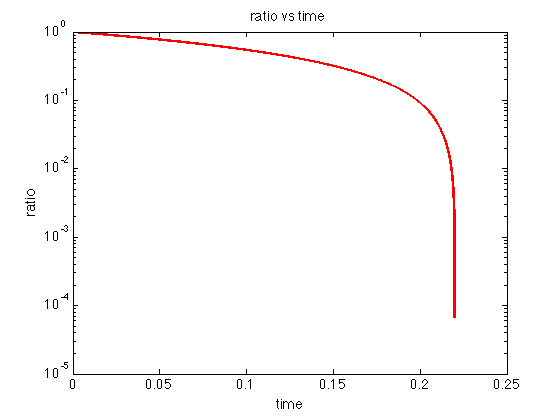}
\includegraphics[width=7.5cm, height=6.5cm]{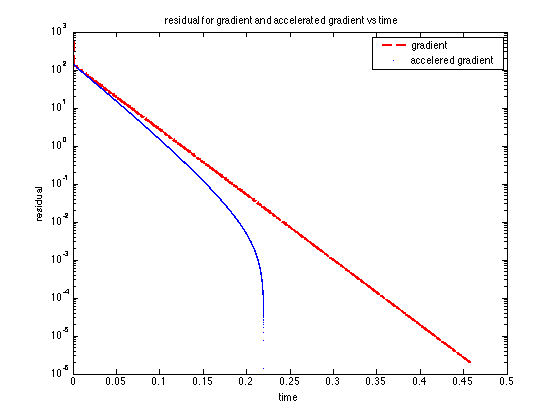}\\
\caption{Acceleration of the convergence in time for the Heat equation with steepest descent method  $n=127$, $X_0=0$, $b=1$, $M$ is built as an incomplete $LU$ factorization of $A$}
\label{ACC_SYSLIN_GRAD2}
\end{center}
\end{figure}

\noindent {\bf Example 5:} Computing the inverse of a matrix\\
We now consider the computation of the inverse of a matrix $A$. We first describe the formal procedure in the general case, we will consider sparse approximations in a second time.
The ODE to be considered here is
\begin{equation}\begin{split}
\Frac{dX}{dt}=Id-AX,\\
X(0)=X_0,
\end{split}\end{equation}
for $A$, e.g. SPD, we have $X(t)\rightarrow A^{-1}$ as $t \rightarrow + \infty$.
\begin{center}
\begin{minipage}[H]{12cm}
  \begin{algorithm}[H]
    \caption{Acceleration for  inverse matrix approximation}\label{AIMA}
    \begin{algorithmic}[1]
    \State $X^{(0)}=X_0$ given in ${\cal M}_n({\mathbb R})$
        \State  ${\tilde X} \approx X^*$ is approached by solving $M{\tilde X} =Id$
        \State $e={\tilde X}-X^{(0)}$ and $u^{(0)}=e$
            \For{$k=0,1, \cdots$ until convergence}
             \State Set $X^{(k+1)}=X^{(k)} +\Delta t (Id-AX^{(k)} +u^{(k)})$
              \State Set   $u^{(k+1)}=(Id-\Delta t A)u^{(k)}$ 
            \EndFor
    \end{algorithmic}
    \end{algorithm}
\end{minipage}
\end{center}
We illustrate here the acceleration obtained by using the accelerator $U$ in two situations. In
Figure \ref{ACC_MAT_GRAD1}, we start from an approximation of the initial error
$e\approx A^{-1}-X_0$, obtained by a threesholding of the coefficients of $A^{-1}$ at a level of $1$ percent; $A$ being the Laplacian finite differences matrix on a square, rescaled by its Frobenius norm ($A=\Frac{A}{\parallel A\parallel_F})$.
We observe clearly the acceleration property of the method.
In Figure \ref{ACC_MAT_GRAD2}, we started from the exact initial error $e=A^{-1} -X_0$ and we observed a much important acceleration effect, this traduces a sensibility. Here $n$ is the number of discretization points in each direction of the domain, the matrix $A$ is then of size $n^2 \times n^2$.
\begin{figure}[!h]
\begin{center}
\includegraphics[width=7.5cm, height=6.5cm]{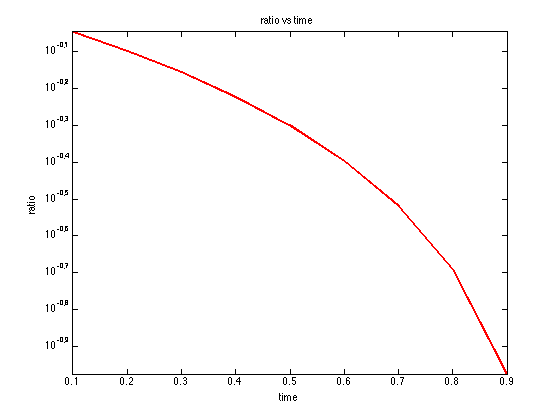}
\includegraphics[width=7.5cm, height=6.5cm]{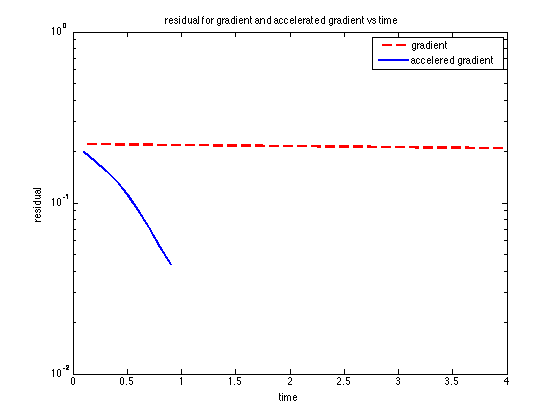}\\
\caption{Acceleration of the convergence in time for inverse of a matrix  $n=63$, $X_0=0, b=1$,$\Delta t=0.1$ and
$e\approx A^{-1}-X_0$}
\label{ACC_MAT_GRAD1}
\end{center}
\end{figure}
\begin{figure}[!h]
\begin{center}
\includegraphics[width=7.5cm, height=6.5cm]{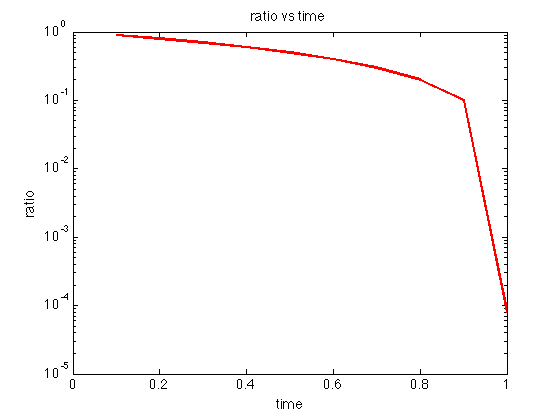}
\includegraphics[width=7.5cm, height=6.5cm]{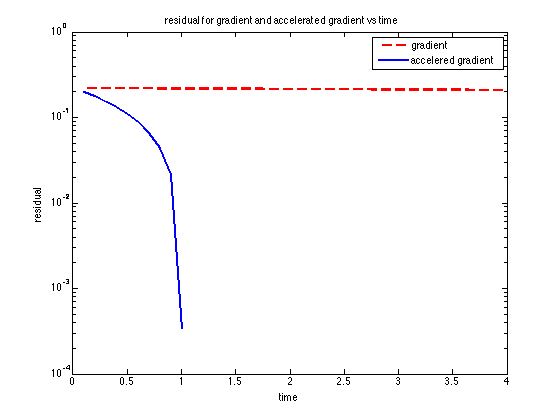}\\
\caption{Acceleration of the convergence in time for inverse of a matrix  $n=63$, $X_0=0, b=1$, $\Delta t=0.1$,
and $e= A^{-1}-X_0$}
\label{ACC_MAT_GRAD2}
\end{center}
\end{figure}


\section{Application of optimal control techniques for capturing steady states}\label{control}
\subsection{Derivation of the method}
A particular acceleration procedure applied to (\ref{ODEaccconv}) consists in choosing $V(t)$ such that $f(A)$ is reached (in a sense to precise) at a given finite time $T$; $V(t)$ appears as a control, defining $X_V(t)$ as the solution of
 \begin{equation}\begin{split}
& \Frac{d X_V}{dt}={\mathcal F}(X_V,A) + V(t) , \ t \in (0,1),\\
 & X_V(0)=X_0.
\end{split}\end{equation}
Hence, $X_V(1)=f(A)$ or ${\mathcal F}(X_V(1),A)=0$. This last condition can be relaxed in terms of optimal control, choosing the merit function
$$
J(X_V,V)=\alpha \parallel {\mathcal F}(X_V(1),A)\parallel^2_F +\gamma \int_0^1 \parallel V(t)\parallel^2_Fdt.
$$
The problem can be rewritten as
$$
\displaystyle{\inf_{V / \Frac{d X_V}{dt}={\mathcal F}(X_V,A) + V(t)} J(X_V,V)}.
$$
In \cite{maday-turinici},  Maday and Turinici proposed a parareal method to solve control problem for parabolic equations; in practice, they used the PM as a preconditioner of the propagation matrix on the fine grid in time; 
in  a recent article \cite{DuSarkisSchaererSzyld}, Du {\it et al} applied the PM procedure as a preconditioning of the Schur complement derived from KKT's optimality relations, for the same type of problem.\\

Let us follow \cite{maday-turinici}: we write for simplicity  $V(t)=Bu(t)$ where $B$ a matrix given in ${\cal M}_{n,p}$. It is proposed to rewrite the problem in terms of virtual control by adding a penalty term to the merit function
$$
J_{\epsilon}(u,\Lambda)=\Frac{\alpha}{2} \parallel  {\mathcal F}(y_{N-1}(T),A)\parallel^2_F +\Frac{1}{2} \int_0^1 \parallel u(t)\parallel^2_Fdt
+\Frac{1}{2\epsilon \Delta t}\displaystyle{\sum_{k=0}^{N-1} \parallel y_{k-1}(T_k^{-})-\lambda_k\parallel^2},
$$
where $\Lambda=(\lambda_0,\cdots, \lambda_{N-1})$ and $\epsilon >0$ is small and $y_k$ is the solution of
\begin{equation}\begin{split}
& \Frac{dy_k}{dt}={\mathcal F}(y_k,A) + Bu(t) , \ t \in (T_k,T_{k+1}),\\
& y_k(T^{-}_k)=\lambda_k,
 \end{split}\end{equation}
with $T_k=k\Delta t$.
The resolution of the minimization $\inf_{u} J_{\epsilon}(u,\Lambda)$ problem can be done by applying, e.g., a gradient method to $J_{\epsilon}(u,\Lambda)$ as follows
$$
\begin{array}{ll}
u_k^{m+1}&=u_k^m-\rho(u_k^m+B^*p_k^m),\\
\lambda_k^{m+1}&=\lambda_k^m-\rho [{\tilde M}^{-1}({\tilde M}^{-1})^*](p_k^m(T_k^+)-p_{k-1}(T^{-}_k))),
\end{array}
$$
where $\rho >0$ is the descent parameter (a conjugate gradient method could be also considered) and where $p_k$ are the solutions of the adjoint equations
\begin{equation}\begin{split}
& -\Frac{dp_{N-1}}{dt}={\mathcal F}^*(p_{N-1},A)  , \ t \in (T_k,T_{k+1}),\\
& p_k(T^{-}_{k+1})=\alpha {\mathcal F}(y_{N-1},A),
\end{split}\end{equation}
and for $k=N-2,N-3,\cdots, 0$
\begin{equation}\begin{split}
& -\Frac{dp_k}{dt}={\mathcal F}^*(p_k,A)  , \ t \in (T_k,T_{k+1}),\\
& p_k(T^{-}_{k+1})=\Frac{1}{\epsilon\Delta t}(y_k(T^{-}_{k+1})-\lambda_{k+1});
\end{split}\end{equation}
 $M$ and ${\tilde M}$ are respectively the fine and the coarse propagators of the PM method:
$$
M=\left(
\begin{array}{llll}
Id & 0 & \cdots & 0\\
-{\bf F} & Id & 0 &\cdots\\
0&-{\bf F} & Id & \cdots \\
0&\cdots&-{\bf F}& Id\\
\end{array}
\right)
\mbox{ and }
{\tilde M}=\left(
\begin{array}{llll}
Id & 0 & \cdots & 0\\
-{\bf G} & Id & 0 &\cdots\\
0&-{\bf G} & Id & \cdots \\
0&\cdots&-{\bf G} & Id\\
\end{array}
\right).
$$
Here ${\bf F}$ and ${\bf G}$ represent the fine and the coarse numerical propagators, on the time intervals $[T_k,T_{k+1}]$, $k=0,\cdots N-1$, which we have not specified since there is no ambiguity.
\subsection{Sparse case}
With the same framework, it is also possible to compute sparse approximations. To this end, we add an additional condition on the matrix coefficients.
Let 
${\cal H}=\{ (i,j) \in I\times J\}$, where $I$ and $J$ are subset on indicies of $\{1, \dots, n\}$. The problem reads as
$$
\displaystyle{\inf_{V / \Frac{d X_V}{dt}={\mathcal F}(X_V,A) + V(t) , \ X_V\in {\cal H}} J(X_V,V)}.
$$

\noindent {\bf Acknowledgements:} This work was supported by a project LEA (Laboratoire Europ\'een Associ\'e CNRS Franco-Roumain). The authors thank Martin J. Gander for fruitful discussions.

E-mail address: {\tt Jean-Paul.Chehab@u-picardie.fr}\\
E-mail address: {\tt Madalina.Petcu@math.univ-poitiers.fr}.

\begin{thebibliography}{99}
\bibitem{BenziRazouk} M. Benzi and N. Razouk, Decay Bounds and O(N) Algorithms for Approximating Functions of Sparse Matrices, ETNA (Electronic Transactions on Numerical Analysis), 28 (2007), pp. 16--39. Special volume in honor of Gene Golub.
\bibitem{bouyouli}
{ R.~Bouyouli, K.~Jbilou, R.~Sadaka and H.~Sadok}, 
{\ Convergence properties of some block Krylov subspace methods for multiple linear systems},
J. Comput. Appl. Math., 196 (2006), pp.~498--511.
\bibitem{jpc} J.-P. Chehab, Differential equations and inverse preconditioners, Computational and Applied Mathematics, Vol 26, N1, pp. 1--34, (2007).
\bibitem{cl} J.-P. Chehab and J. Laminie, Differential equations and solution of linear systems, Numerical Algorithms, 40, pp.~103--124 (2005).
\bibitem{AbreuRaydan} B. De Abreu and M. Raydan, Residual methods for the large-scale matrix p-th root and some related problems, Applied Mathematics and Computation, Vol. 217, pp.~650--660, (2010).
\bibitem{DuboisSaidi} F. Dubois and A. Saidi, Unconditionally stable scheme for Riccati equation, ESAIM Proc,
8 (2000), pp.~39--52.
 \bibitem{DuSarkisSchaererSzyld} X. Du, M. Sarkis, C. E. Schaerer and D. B. Szyld, Inexact and truncated Parareal-in-time Krylov subspace methods for parabolic optimal control problems, Electronic Transactions on Numerical Analysis, vol. 40 (2013), pp. 36--57.
\bibitem{farhat-03}
{ C.~Farhat and M.~Chandesris}, { Time-decomposed parallel
time-integrators: theory and feasibility studies for fluid, structure,
and fluid-structure applications}, Internat. J. Numer. Methods Engrg.,
58 (2003), pp.~1397--1434.
\bibitem{farhat-06}
{ C.~Farhat, J.~Cortial, C.~Dastillung and H.~Bavestrello}, {
Time-parallel implicit integrators for the near-real-time prediction
of linear structural dynamic responses}, Internat. J. Numer. Methods
Engrg., 67 (2006), pp.~697--724.
\bibitem{Fommer} A. Frommer, S. Guettel and M. Schweitzer,Efficient and stable Arnoldi restarts for matrix functions based on quadrature, SIAM J. Matrix Anal. Appl., 35:661--683, 2014
\bibitem{GallopoulosSaad} E. Gallopoulos and Y. Saad, On the parallel solution of parabolic equations, Proc. 1989
ACM Internat. Conf. on Supercomputing, Heraklion, Greece, 1989, pp. 17--28.
\bibitem{GanderGuettelParaexp} M. Gander and S. Guettel, PARAEXP: A parallel integrator for linear initial-value problems, SIAM J. Sci. Comput., 35(2):C123--C142, 2013
\bibitem{GanderPetcu1} M. Gander and M. Petcu, Analysis of the modified parareal algorithm for second-order ordinary differential equations, AIP Conference Proceedings: Numerical Analysis and Applied Mathematics, pp.~233--236, (2007).
\bibitem{GanderPetcu2} M. Gander and M. Petcu, Analysis of a Krylov subspace enhanced parareal algorithm for linear problems. ESAIM Proc, 25, pp.~114--129, November 25, (2008).
\bibitem{GanderVdW}M. Gander and S. Vandewalle, Analysis of the parareal time-parallel time-integration method, SIAM J. Sci. Comput.,
Vol. 29, No. 2, pp. 556--578, (2007).
\bibitem{GuettelRational} S. Guettel, Rational Krylov approximation of matrix functions: Numerical methods and optimal pole selection, GAMM Mitteilungen, 36(1):8--31, 2013.
 \bibitem{HighamSine} G. I. Hargreaves and N. J. Higham , Efficient algorithms for the matrix cosine and sine, Numerical Algorithms (2005) 40, pp.~383--400.
 \bibitem{HighamBook} N. Higham, Functions of Matrices: Theory and Computation,
 SIAM, 2008. xx+425 pages.
 \bibitem{Helmke}  U. Helmke and J.B. Moore, Optimization and Dynamical Systems, Comm. Control Eng.
Series, Springer, London, 1994.
 \bibitem{Jbilou1} K. Jbilou, A. Messaoudi and H. Sadok, Global FOM and GMRES algorithms for matrix equations, Applied Numerical Mathematics 31 (1999), pp.~49--63
\bibitem{lions-maday-turinici}
{J.-L. Lions, Y.~Maday and G.~Turinici}, { R\'esolution
d'{EDP} par un sch\'ema en temps ``parar\'eel''},
C. R. Acad. Sci. Paris S\'er. I Math., 332 (2001), pp.~661--668.
\bibitem{maday-turinici}
{ Y.~Maday and G.~Turinici}, { Parallel in time algorithms for
quantum control: the parareal time discretization scheme},
Int. J. Quant. Chem., 93 (2003), pp.~223--228.
\bibitem{MadayTuriniciCRAS}  Y.~Maday and G.~Turinici,  A parareal in time procedure for the control of partial differential equations. C. R. Math. Acad. Sci. Paris 335 (2002), no. 4, pp.~387--392.
 \bibitem{MolerExponential} C.~Moler and C.~V.~Loan, Nineteen Dubious Ways to Compute the Exponential of a Matrix, Twenty-Five Years Later, SIAM REVIEW, (2003)  Vol. 45, No. 1.
 \bibitem{MonsalvesRaydan1} M Monsalve and M. Raydan, A new inversion-free method for a rational matrix equation
 Linear Algebra and Its Applications, Vol. 433, no 1, pp 64--71, 2010
 \end{thebibliography}
\end{document}